\documentclass[a4paper,11pt]{article}
\usepackage{geometry,amssymb,amsmath,enumerate,latexsym,array}
\usepackage[dvips]{graphicx}
\usepackage[all]{xy}
\begin{document}
\newtheorem{example}{Example}[section]
\newtheorem{alg}[example]{Algorithm} 
\newtheorem{rem}[example]{Remark} 
\newtheorem{lem}[example]{Lemma} 
\newtheorem{thm}[example]{Theorem} 
\newtheorem{cor}[example]{Corollary} 
\newtheorem{Def}[example]{Definition} 
\newtheorem{prop}[example]{Proposition} 
\newenvironment{proof}{\noindent {\bf Proof} }{\mbox{} \hfill $\Box$
\mbox{}\\}

 \title { One-sided Noncommutative Gr\"obner Bases with\\
           Applications to Computing Green's Relations
 \thanks{KEYWORDS: Gr\"obner bases, semigroups, right ideals, left ideals. 
 \newline AMS 1991 CLASSIFICATION: 16-04 16D15  }}
\author{ Anne Heyworth \thanks{ Supported 1995-8 by an EPSRC
Earmarked Research Studentship, `Identities among relations for
monoids and categories', and 1998-9 by  a University of Wales, Bangor,
Research Assistantship.}\\ School of Mathematics \\ University of Wales,
Bangor \\
Gwynedd LL57 1UT \\ United Kingdom \\ map130@bangor.ac.uk} 
\maketitle

\begin{center} Abstract \end{center}
{\small 
Standard noncommutative Gr\"obner basis procedures are used for computing
ideals of free noncommutative polynomial rings over fields.
This paper describes Gr\"obner basis 
procedures for one-sided ideals in finitely presented noncommutative 
algebras over fields. 
The polynomials defining a $K$-algebra $A$ as a quotient 
of a free $K$-algebra are combined with the polynomials defining a one-sided 
ideal $I$ of $A$, by using a tagging notation. 
Standard noncommutative Gr\"obner basis techniques can then be applied to
the mixed set of polynomials, thus calculating $A/I$ whilst working in a free
structure, avoiding the complication of computing in $A$.
The paper concludes by showing how the results can be applied to
completable presentations of semigroups and so enable calculations of 
Green's relations.}

\section{Introduction}

In 1965 Buchberger invented Gr\"obner basis theory, techniques enabling the 
computation of ideals in commutative polynomial rings over fields.
Implementations of Buchberger's algorithm are now provided by all major
computer algebra systems, a good cross-section of the ways in which the 
theory has developed may be found in \cite{RISC}. 
Mora generalised Gr\"obner basis theory to noncommutative polynomial rings
(algebras) \cite{FMora}. 
Introductions to these procedures may be found in \cite{Ufn,TMora2}.
This paper presents an extension of the noncommutative Gr\"obner
basis procedures for polynomials to what we call {\em tagged polynomials}.
The intention is to describe methods of computation that may be applied
to the problem of computing right (or left) ideals in finitely presented
$K$-algebras.
\\

The data defining the problem consists of the field $K$, a set of
noncommuting
variables $X$, a set of generators $P \subseteq K[X^\dagger]$ for a two-
sided
ideal $\langle P \rangle$, defining an algebra 
$A=K[X^\dagger]/ \langle P \rangle$ and a set of generators $Q'
\subseteq A$
for a right ideal $\langle Q' \rangle^r$.
We expect elements of $A$ to be given in terms of $K[X^\dagger]$, so
$Q'$ is
specified by a set $Q \subseteq K[X^\dagger]$.
The problem we address is that of computing the right $A$-ideal
generated by $Q'$, written $\langle Q' \rangle^r$.
\\

Our solution lies in using the free right $K$-module $K[\dashv \!
X^\dagger]$. 
Here  $\dashv$ is just a symbol or tag and 
$K[\dashv\!X^\dagger]$ is bijective with $K[X^\dagger]$.
We call elements of $K[\dashv\!X^\dagger]$ \emph{tagged polynomials} 
(Definition 3.1) and
write them $k_1\dashv\!m_1+ \cdots + k_n\dashv\!m_n$ where 
$k_1,\ldots,k_t \in K$ and $m_1,\ldots,m_t \in X^\dagger$. 
Ordinary polynomials $F_P$ defining the two-sided ideal
$\langle P \rangle$ which defines $A$ are combined with
tagged polynomials $F_T$ defining the one-sided ideal. 
The mixed set of polynomials $F:=(F_T,F_P)$ determines a
reduction relation $\to_F$ (Definition 3.2) on the tagged polynomials
$K[\dashv \! X^\dagger]$.
\\

The value of this combination and use of tagging is in computation, as
will be shown in 
the main result (Algorithm 4.9), which describes a variant of the Buchberger
algorithm. The initial mixed set of polynomials $F$ is appended with tagged 
and non-tagged polynomials until the relation $\to_F$ is complete on
$K[\dashv\!X^\dagger]$. When the procedure terminates the usual
normal form arguments apply and reduction modulo $F$ can be used to
solve 
the membership problem for the right ideal $\langle Q' \rangle^r$
of the finitely presented algebra $A$.
\\
 
Previous work \cite{Birgit,MRC} attempt the computation of one-sided ideals
by using different definitions of 
purely one-sided reduction relation in particular algebras 
(e.g. $\mathbb{Q}[M]$ for a monoid $M$ presented by a complete rewrite system). 
The main problem encountered is that of computing in a non-free algebra,
we avoid this and base all the computations specifying the algebra at
the same level (in a particular free right module) 
as those for the ideal and compute the two simultaneously. 
In other words, the methods we describe 
provide for local computations, concerning single ideals 
$\langle Q' \rangle^r$ without the requirement to compute 
the global structure of the algebra $A$ or face the
difficulties of calculations with elements of $A$.
This idea follows the philosophy that 
computations take place in free objects.

\section{Algebra Presentations and One-sided Ideals}

If $X$ is a set, then 
$X^\dagger$ is the {\em free semigroup} 
of all strings of elements of $X$, and
$X^*$ is the {\em free monoid} of all strings 
together with the empty string, which acts as the identity $i\!d$ for
$X^*$.
A {\em semigroup presentation} is a pair $sgp \langle X | R \rangle$
where
$X$ is a set and $R \subseteq X^\dagger \times X^\dagger$. 
It {\em presents a semigroup} $S$ if $X$ is a set of generators of $S$
and the natural morphism $\theta:X^\dagger \to S$ induces an isomorphism
from $X^\dagger\!/=_R$ to $S$, where $=_R$ is the
congruence generated on $X^\dagger$ by $R$.
Similarly, a {\em monoid presentation} is a pair $mon \langle X | R \rangle$
where $X$ is a set and $R \subseteq X^* \times X^*$. 
It {\em presents a monoid} $M$ if $X$ is a set of generators of $M$
and the natural morphism $\theta:X^* \to M$ induces an isomorphism
from $X^*\!/=_R$ to $M$, where $=_R$ is the
congruence generated on $X^*$ by $R$.
\\

Let $K$ be a field.
The { \em free $K$-algebra} $K[S]$ on a semigroup $S$ 
consists of all sums of $K$-multiples
of elements of $S$ with the operations of 
addition and multiplication defined
in the obvious way. In particular the elements of
$K[X^\dagger]$ are called {\em polynomials} and written
$k_1m_1+\cdots+k_nm_n$ where $k_1,\ldots,k_n \in K$ and 
$m_1,\ldots,m_n \in X^\dagger$.\\

If $P$ is a subset of an algebra $Z$ then the 
{\em two-sided ideal} generated by $P$ in 
$Z$ is denoted $\langle P \rangle$. In the case $Z=K[X^\dagger]$ this
consists of all sums of multiples of
elements of $P$, i.e.
$$\langle P \rangle := \{ k_1u_1p_1v_1+ \cdots + k_nu_np_nv_n 
               \ | \ p_1, \ldots,p_n \in P, k_1,\ldots,k_n \in K, 
               u_1,v_1,\ldots,u_n,v_n \in X^*\}.
$$
Given an ideal in an algebra the {\em membership problem} is that of
determining,
for a given element of the algebra, whether it is an element of the
ideal.
\\

A {\em congruence} on an algebra $Z$ is an equivalence relation $\sim$
on its elements such that if $p \sim q$ then $p+u \sim q+u$ and 
$upv \sim uqv$ for all $u,v \in Z$.
Given an algebra $Z$ and an ideal $\langle P \rangle$ 
ideal membership defines a congruence on the algebra by
$$
p \sim q \Leftrightarrow p-q \in \langle P \rangle.
$$
The {\em quotient algebra} $Z/\langle P \rangle$ is the algebra of
congruence 
classes of $Z$ under $\langle P \rangle$.
A $K$-\emph{algebra presentation} is a pair $alg\langle X | P \rangle$, where
$P \subseteq K[X^\dagger]$.
A $K$-algebra $A$ is presented by $alg\langle X|P \rangle$ if 
$X$ is a set of generators of $A$ and the natural morphism
$K[X^\dagger] \to A$ induces an isomorphism $K[X^\dagger]/\langle P
\rangle \to A$.\\

Noncommutative Gr\"obner basis theory (as described in \cite{TMora, TMora2})
uses the notion of an ordering on $X^\dagger$, thereby allowing the 
concepts of {\em leading monomial}, 
{\em leading term} and {\em leading coefficient} 
 on the polynomials of $K[X^\dagger]$. 
Given any subset $P$ of $K[X^\dagger]$ 
an ordering determines a Noetherian 
reduction relation $\to_P$ on the
elements of $K[X^\dagger]$. 
The reflexive, symmetric, transitive closure of this
relation is a congruence relation 
coinciding with the ideal membership of $\langle P \rangle$.
\\

Let $A$ be the $K$-algebra presented by
$alg\langle X | P \rangle$ and let $Q' \subseteq A$.
We wish to consider  the right ideal $\langle Q'\rangle^r $
generated in $A$ by $Q'$,
i.e.
$$
\langle Q'\rangle^r :=
  \{q_1'a_1+ \cdots +q_n'a_n \ | \ a_1, \ldots, a_n \in A,
                                           q_1', \ldots, q_n' \in Q'\}.
$$

A {\em right congruence} on an algebra $A$ is an equivalence relation
$\stackrel{r}{\sim}$ such that for all $a,b,y \in A$
$$a \stackrel{r}{\sim} b  \Rightarrow 
       a+y \stackrel{r}{\sim} b+y  \text{ and } 
       ay \stackrel{r}{\sim} by.$$

Membership of a right ideal $\langle Q' \rangle^r$ defines 
a right congruence on $A$,  
by $a \stackrel{r}{\sim}_{Q'} b \Leftrightarrow a-b \in \langle Q'
\rangle^r$.
The quotient $A/\langle Q' \rangle^r$ is 
the set of all the right congruence classes of $A$ under 
$\stackrel{r}{\sim}_{Q'}$ 
where classes are denoted $[a]_{Q'}$ for $a \in A$.
Note that for $a, b \in A$,
$[a+b]_{Q'} = [a]_{Q'} + [b]_{Q'}$ 
and
$[a]_{Q'}[b]_{Q'} = [a]_{Q'}$.\\

 {\em Buchberger's algorithm} 
is a critical pair completion procedure. 
The algorithm begins with a set of polynomials 
$P$ of a free algebra. Set $F:=P$ and a search for 
overlapping leading terms 
will find all critical terms of the 
reduction relation $\to_F$. 
This enables a test for local confluence. 
Overlaps which cannot be resolved result in 
{\em S-polynomials} all of which are 
added to $F$ at each stage
(though some elimination is possible, 
for efficiency). The algorithm
terminates when all the overlaps of 
$F$ can be resolved, i.e.
$\to_F$ is complete 
(Noetherian and confluent), when
this occurs $F$ is said to be a 
{\em Gr\"obner basis} for 
the ideal $\langle P \rangle$.
Obtaining a Gr\"obner basis $F$ 
allows, in particular, the 
solution of the  membership problem 
by using $\to_F$ as a normal
form function on $K[X^\dagger]$.
Thus, if $F$ is a Gr\"obner basis for 
the ideal $\langle P \rangle$ on $K[X^\dagger]$ 
and $p,q \in K[X^\dagger]$, then
$$
  p \sim q \Leftrightarrow 
  p-q \in \langle P \rangle \Leftrightarrow
p \stackrel{*}{\to}_{F} u \text{ and }
q \stackrel{*}{\to}_{F} u 
  \text{ for some } u \in K[X^\dagger].
$$

The Noetherian property ensures that 
the process of reduction terminates with an 
irreducible element; confluence ensures that 
any two elements of the same class
reduce to the same form. In practice: reduce 
each polynomial as far as 
possible using $\to_{F}$, the original 
polynomials are equivalent only 
if their irreducible forms are equal.\\

In the next sections we show how to apply Buchberger's algorithm to obtain --
when possible -- a Gr\"obner basis of (two types of) polynomials, which will
enable the use of normal forms arguments.

\section{One-sided Noncommutative Gr\"obner Basis Procedures}

Given a finitely presented $K$-algebra $A$ and a subset $Q'$ of $A$ we
wish to compute the right ideal $\langle Q' \rangle^r$. 
The meaning of `computing the ideal' in this context is that of solving  
the ideal membership problem for $\langle Q' \rangle^r$ in $A$.
The $K$-algebra $A$ is presented by 
$alg\langle X | P \rangle$ and to obtain normal forms for $A$ we would
therefore apply Gr\"obner basis procedures
to $P$ in the free algebra $K[X^\dagger]$.
Since we are interested in a one-sided ideal we introduce the 
tagging notation which will allow the combination of $P$ and $Q$.

\begin{Def}[Tagged polynomials]
Let $K$ be a field, let $X$ be a set and let $\dashv$ be a symbol.
Then $\dashv\!X^\dagger$ is the set of \textbf{tagged terms}
$\dashv\!m$
where $m \in X^\dagger$ and $K[\dashv \! X^\dagger]$ is the free right 
$K[X^\dagger]$-module of \textbf{tagged polynomials},
i.e. elements $k_1 t_1 + \cdots + k_n t_n$ for
$k_1, \ldots, k_n \in K$, $t_1, \ldots, t_n \in \dashv\! X^\dagger$.
\end{Def}

Let $>$ be a {\em semigroup ordering} on $X^\dagger$, i.e.
$>$ is an irreflexive, antisymmetric, transitive relation on $X^\dagger$
such
that if $m_1>m_2$ then $um_1v>um_2v$ for all $u,v \in X^*$.
Further we require the {\em well-ordering} property, that there is no
infinite sequence $m_1>m_2>m_3> \cdots$.\\

Let $p=k_1m_1 + \cdots + k_nm_n \in K[X^\dagger]$. 
The $k_im_i$ are referred to as the {\em monomials} of the polynomial, where
$m_i$ is the {\em term} and $k_i$ the {\em coefficient}.
Assuming the well-ordering on $X^\dagger$, 
the leading monomial $\mathtt{LM}(p)$ is defined to be the monomial with
the largest term. The leading term $\mathtt{LT}(p)$ and
leading
coefficient $\mathtt{LC}(p)$ are the coefficient and term of
this monomial.\\ 

To simplify the definitions throughout this paper 
we will assume all polynomials to be {\em monic},
i.e. their leading coefficients are all 1.
There is no loss in doing this: $K$ is a field so 
the polynomials $F$ may always be divided by their leading 
coefficients and still generate the same ideal.\\

The well-ordering on $X^\dagger$ induces a well-ordering on
$\dashv\!X^\dagger$
defined by $\dashv\!m_1>\dashv\!m_2 \Leftrightarrow m_1>m_2$.
This gives corresponding notions of leading monomial, leading term 
and leading coefficient for the tagged polynomials. 
In detail: if $p=k_1m_1+ \cdots + k_nm_n$ where $k_1,\ldots,k_n \in K$ and
$m_1,\ldots,m_n\in X^\dagger$ is a polynomial with leading term
$\mathtt{LT}(p)=m_i$ then the tagged polynomial
$\dashv\!p:=k_1\dashv\!m_1 + \cdots + k_n \dashv\!m_n$ has a tagged leading 
term $\mathtt{LT}(\dashv\!p)=\dashv\!m_i$. 
\\

We will now introduce the definition of a reduction relation 
on $K[\dashv\!X^\dagger]$, defined by a mixed set of polynomials
$F=(F_T,F_P)$ where $F_T$ is a set of tagged polynomials, elements of
the module, and $F_P$ is a set of polynomials, 
elements of the algebra acting on the right of the module.
The reduction relation $\to_F$ combines the two relations so
that they are defined on the free right module of tagged polynomials.

\begin{Def}[Reduction of tagged polynomials]
Let $F:=(F_T,F_P)$ where 
$F_T \subseteq K[\dashv\! X^\dagger]$ and 
$F_P \subseteq K[X^\dagger]$. 
Define the reduction relation $\to_F$ on
tagged polynomials $f \in K[\dashv \! X^\dagger]$ by
$$
f \to_{F}  f - k u (f_i) v
$$
if $u \mathtt{LT}(f_i) v$ occurs in $f$ with coefficient $k \in K$
for some $u \in \dashv\!X^* \cup \{i\!d\}$, $v \in X^*$, $f_i \in F$.  
\end{Def}

A one-step reduction like that of the definition may also be
written $f \to_{f_i}  f - k u (f_i) v$. 
This relation may be understood to be a rewrite
system on the polynomials (similarly to observations made in \cite{Birgit}
on Mora's definitions of reduction \cite{FMora}). When a multiple of the 
leading term of $f_i$ for $f_i \in F$ occurs in the polynomial that is 
to be reduced, the rest of $f_i$ is substituted for the leading monomial 
of $f_i$.

Regarding $F$ as a rewrite system with two types of rules that may be
applied to monomials of polynomials, we could say that
the non-tagged polynomials can be applied anywhere in a term, 
but the tagged ones apply only at the tagged side of a term.

\begin{example}[Reduction]
\emph{For example let $F_T:=\{f_1,f_2\}$ where 
$f_1:=\dashv\!xyx+\dashv\!yx+2\dashv\!y$,
$f_2:=\dashv\!yx^2+\dashv\!x^2$ and
$F_P:=\{f_3,f_4\}$ where $f_3:=x^2y-3yx$, $f_4:=yx^3-2xy$. 
Then the tagged polynomial 
$f:=8\dashv\!xyx^2y^3+5\dashv\!y$ cannot be reduced by $f_2$ or $f_4$ 
but can be reduced by $f_1$ to 
$f-8f_1xy^3=5\dashv\!y-8\dashv\!yx^2y^3-16\dashv\!yxy^3$ 
or by $f_3$ to
$f-8\dashv\!xyf_3y^2=5\dashv\!y+24\dashv\!xy^2xy^2$.}
\end{example}

These results allow the combination of two-sided and one-sided
congruences,
by the use of tagged polynomials. A mixed set of polynomials $F$
defines a reduction relation on the module of all tagged polynomials.
The reflexive, symmetric, transitive closure of $\to_F$ will be 
denoted $\stackrel{*}{\leftrightarrow}_F$. 
The class of $f \in K[\dashv \! X^\dagger]$ under the equivalence
relation 
$\stackrel{*}{\leftrightarrow}_F$ will be denoted $[f]_F$.

\begin{thm}
Let $A$ be a $K$-algebra finitely presented by $alg\langle X | P
\rangle$ with quotient morphism $\theta$.
Let $Q \subseteq K[X^\dagger]$ and define $Q':=\theta Q$.
Define $F:=(\dashv\!Q, P)$ where $\dashv\!Q:=\{\dashv\!q:q\in Q\}$.
Then there is a bijection of sets
$$
\frac{ K[\dashv \! X^\dagger] }{ \stackrel{*}{\leftrightarrow}_F }
\cong  \frac{A}{ \langle Q' \rangle^r }$$
\end{thm}

\begin{proof}
The quotient
morphism $\theta:K[X^\dagger] \to A$, 
defines a surjection
$\theta^\dashv:K[\dashv \! X^\dagger] \to A$. 
Then $\theta^\dashv(\dashv \! Q)=Q'$.\\

Define $\phi:K[X^\dagger]/\!\stackrel{*}{\leftrightarrow}_F \; 
            \to A/\langle Q' \rangle^r$ by 
\, $\phi([f]_F) := [\theta^\dashv(f)]_{Q'}.$

To prove that $\phi$ is well-defined we show that it preserves the right 
congruence $\stackrel{*}{\leftrightarrow}_F$.
We assume all polynomials of $F$ are monic.
Let $f \in K[\dashv\!X^\dagger]$
and $f_i \in F$ and suppose that $f \to_F f - kuf_iv$ for some
$k \in K$, $u \in \dashv\!X^* \cup \{i\!d\}$, $v \in X^*$.
By definition 
$\phi([ f-kuf_iv ]_F) = [ \theta^\dashv(f) - \theta^\dashv(kuf_iv)
]_{Q'}$.
\\

Now either\\ 
(i) $f_i \in P \subseteq K[X^\dagger]$ and  
      $\theta^\dashv(kuf_iv)=0$, since $kuf_iv \in \langle P \rangle$,\\
or else \\
(ii) $f_i \in \dashv\! Q \subseteq K[\dashv\!X^\dagger]$ and 
      $u=i\!d$, so $\theta^\dashv(kuf_iv)=k\theta^\dashv (f_iv)$.\\
In either case $\theta^\dashv(kuf_iv) \in \langle Q' \rangle^r$, so
$\phi(  [f ]_F ) = \phi( [ f-kuf_iv ]_F )$, 
i.e. $\phi$ preserves the relation $\to_F$. 

Furthermore if $[f]_F=[g]_F$ for some $f,g\in K[\dashv\!X^\dagger]$, then
for all $v \in X^*$,
\begin{align*} 
\phi([fv]_F) &= [\theta^\dashv(fv)]_{Q'}\\
             &= [\theta^\dashv(f)]_{Q'} \theta(v) \\
             &= [\theta^\dashv(g)]_{Q'} \theta(v) \\
             &= [\theta^\dashv(gv)]_{Q'} \\
             &= \phi([gv]_F).
\end{align*}
Therefore $\phi$ preserves the right congruence 
$\stackrel{*}{\leftrightarrow}_F$.\\

We now prove that $\phi$ is surjective.
Let $a \in A$. 
Then there exists $f \in K[\dashv \! X^\dagger]$ such that 
$\theta^\dashv(f)=a$, because $\theta^\dashv$ is a surjection. 
Thus for all $[a]_{Q'} \in A/\langle Q' \rangle^r$
there exists $[f]_F 
\in K[\dashv \! X^\dagger]/\stackrel{*}{\leftrightarrow}_F$
such that $\phi([f]_F)=[\theta^\dashv f]_{Q'}=[a]_{Q'}$.\\

Finally, we prove that $\phi$ is injective.
Let $f,g \in K[\dashv \! X^\dagger]$ such that 
$\phi[f]_F = \phi[g]_F$.
Then $[\theta^\dashv(f)]_{Q'} = [\theta^\dashv(g)]_{Q'}$. 
Therefore there exist $q_1',\ldots,q_n' \in Q'$ and $k_1,\ldots,k_n \in
K$, 
$a_1,\ldots,a_n \in A$, such that
$$\theta^\dashv(f)-\theta^\dashv(g)=k_1q_1'a_1+\cdots + k_nq_n'a_n.$$
For $i=1,\ldots,n$ there exists $q_i \in Q$, $y_i \in X^*$ such that
$\theta^\dashv(q_iy_i)=q_i'a_i$.
Hence 
$$\theta^\dashv(f)-\theta^\dashv(g)=k_1q_1y_1+ \cdots + k_nq_ny_n.$$
Now $\theta^\dashv$ preserves $+$ and therefore
$\theta^\dashv( f - g - k_1q_1y_1 - \cdots - k_nq_ny_n ) = 0$.
By the definition of $\theta$, and so $\theta^\dashv$, therefore
$$f-g-k_1q_1y_1 - \cdots - k_nq_ny_n = l_1u_1p_1v_1 + \cdots +
l_mu_mp_mv_m$$
for some $p_1,\ldots,p_m \in P$, 
$l_1,\ldots,l_m\in K$ and $u_1,\ldots,u_m \in \dashv\!X^* \cup \{i\!d\}$, 
$v_1,\ldots,v_m \in X^*$.
Therefore $f \stackrel{*}{\leftrightarrow}_F g$, from the definition.
Therefore $\phi$ is a well-defined bijection of sets.
\end{proof}

\begin{cor}
Let $S$ be a semigroup with presentation $sgp\langle X | R\rangle$.
Let $P:=\{l-r:(l,r)\in R\}$, $Q \subseteq K[X^\dagger]$.
Define $F:=(P,\dashv\!Q)$. Then there is a bijection of sets
$$\frac{K[S]}{\langle Q' \rangle^r} \cong 
  \frac{K[\dashv\!X^\dagger]}{\stackrel{*}{\leftrightarrow}_F}$$
\end{cor}

Here it is appropriate to observe the link to rewrite systems which is used
in the proof of this corollary, in particular, $alg\langle X|P\rangle$ is
a presentation of $K[S]$ \cite{TMora,Birgit,paper3}. 
This corollary (also see the next result) provides an alternative
approach to that of Reinert and Zecker \cite{MRC} for attempting the computation of ideals in 
$\mathbb{Q}[M]$, where $M$ is a monoid. 
Our computations are based in $\mathbb{Q}[X^*]$ where $X$ is a set of 
generators for $M$, the computations of Reinert and Zecker are made within
$\mathbb{Q}[M]$, (also using a presentation of $M$).

\begin{thm}
 Let $X$ be a set of generators for the terms of a $K$-algebra $A$ and let
$P \subseteq K[X^*]$ such that the natural morphism $\theta:K[X^*] \to A$
induces an isomorphism $K[X^*]/\langle P \rangle \to A$.
Let $Q \subseteq K[X^*]$ and define $Q':=\theta Q$.
Define $F:=(\dashv\!Q, P)$ where $\dashv\!Q:=\{\dashv\!q:q\in Q\}$.
Then there is a bijection of sets
$$
\frac{ K[\dashv \! X^*] }{ \stackrel{*}{\leftrightarrow}_F }
\cong  \frac{A}{ \langle Q' \rangle^r }$$
\end{thm}
\begin{proof}
Define 
$\phi:K[\dashv\!X^*]/\stackrel{*}{\leftrightarrow}_F \; \to A/ \langle Q' \rangle^r$ 
by
$\phi([f]_F):=[\theta^\dashv(f)]_{Q'}$.
The verification that $\phi$ is a well-defined bijection on the congruence 
classes is similar to that detailed in the proof of Theorem 3.4.
\end{proof}

\section{The Noncommutative Buchberger Algorithm for One-sided Ideals}

Recall that $F=(F_T,F_P)$ is a mixed set of polynomials 
$F_T \subseteq K[\dashv\!X^\dagger]$ and
$F_P \subseteq K[X^\dagger]$.
The definition of reduction of a tagged polynomial
$f$ requires that a tagged term $\dashv \! m$ of $f$ is some
multiple
of a leading term from the polynomials $f_i$ of $F$. 
This definition of reduction will 
allow the application of the standard noncommutative
Buchberger algorithm to $F$ to attempt to complete 
$\to_F$.

\begin{Def}[Gr\"obner basis of mixed polynomials]
A set $F=(F_T,F_P)$ where $F_T \subseteq K[\dashv\!X^\dagger]$ and
$F_P \subset K[X^\dagger]$ is a \textbf{Gr\"obner basis} on 
$K[\dashv\!X^\dagger]$ with respect to $>$ if $\to_F$ is complete.
\end{Def}

\begin{lem}[Noetherian property]
Let $F=(F_T,F_P)$ where $F_T \subseteq K[\dashv \! X^\dagger]$ and 
$F_P \subseteq K[X^\dagger]$. 
Let $>$ be a semigroup well-ordering on $X^\dagger$.
Then the reduction relation $\to_{F}$ is Noetherian on $K[\dashv \!
X^\dagger]$.
\end{lem}

\begin{proof}
According to the definition, the process of reduction replaces one monomial
with monomials which are smaller with respect to $>$ 
(since $>$ is a term order on $X^\dagger$).
The existence of an infinite sequence of reductions
$f_1 \to_{F} f_2 \to_{F} \cdots$ of polynomials 
$f_1,f_2, \ldots \in K[\dashv \! X^\dagger]$ would imply the existence
of
an infinite sequence $m_1>m_2> \cdots$ of terms 
$m_1,m_2, \ldots \in X^\dagger$.
Therefore $\to_{F}$ is Noetherian.
\end{proof}

\begin{Def}[Matches and S-polynomials of tagged and non-tagged
polynomials]
\mbox{ }\\
Let $F=(F_T,F_P)$ where 
$F_T \subseteq K[\dashv\!X^\dagger]$ and 
$F_P \subseteq K[X^\dagger]$.
A pair of polynomials $f_1,f_2 \in F$ 
has a \textbf{match} if their
leading terms $m_1,m_2$ coincide. If a pair of polynomials have a match then
an \textbf{S-polynomial} is defined. There are five possible cases:
\begin{center}
\begin{tabular}{|l|ll|l|l|}
\hline
      & & match & S-polynomial & \\
\hline
both $f_1$ and $f_2$ in $F_T$  &    (i)    & $m_1 v = m_2$   & $f_1v-f_2$  & where $v\in X^*$ \\
\hline
$f_1$ in $F_T$ and $f_2$ in $F_P$ & (ii)   & $m_1 v = u m_2$ & $f_1v-uf_2$ &  \\ 
& (iii)                                    & $m_1 = u m_2 v$ & $f_1-uf_2v$ & where $u \in \dashv\!X^*\cup\{i\!d\}, v \in X^*$ \\
\hline
both $f_1$ and $f_2$ in $F_P$ & (iv)       & $u m_1 = m_2 v$ & $f_1v-uf_2$ & \\
&  (v)                                     & $m_1 = u m_2 v$ & $f_1-uf_2v$ & where $u,v \in X^*$\\
\hline
\end{tabular}
\end{center}

A match is said to \textbf{resolve} if the resulting S-polynomial can be
reduced to zero by $F$.
\end{Def}

\begin{rem}
{\em
If a match of any of the types above occurs between $f_1$ and $f_2$
then the match may be represented:
$u_1 m_1 v_1 =u_2 m_2 v_2$, where $u_1,u_2,v_1,v_2 \in \dashv\!X^* \cup X^*$.
A match of $f_1$ and $f_2$ may occur when either, neither, or both of
$f_1$ and $f_2$ are tagged.
However, if one or both has a tag, the tag forms part of the match and
the resulting S-polynomial will be tagged.}
\end{rem}

The following lemma is proved in the same way as in the standard 
commutative non-tagged situation as described in \cite{TAT}.

\begin{lem}
\label{zero}
Let $F=(F_T,F_P)$ where $F_T \subseteq K[\dashv \! X^\dagger]$
and $F_P \subseteq K[X^\dagger]$.
Let $g_1, g_2 \in K[\dashv\!X^\dagger]$ where $g_1-g_2 \to_{F}^* 0$. 
Then there exists a tagged polynomial $h \in K[\dashv\!X^\dagger]$ such
that 
$g_1 \stackrel{*}{\to}_{F} h$ and
$g_2 \stackrel{*}{\to}_{F} h$.
\end{lem}

\begin{proof} 
The {\em length} of a reduction sequence is 
defined to be the number of one-step reductions of which it is made up.
This proof is by induction on the length of the reduction sequence
$g_1-g_2 \stackrel{*}{\to}_{F} 0$.\\

For the basis of induction 
suppose the length of the reduction sequence is zero.
Then $g_1-g_2=0$ so $g_1=g_2$.\\

For the induction step, assume that if $g_1'-g_2' \stackrel{*}{\to}_{F}
0$
is a reduction sequence of length $n$ then there exists 
$h \in K[\dashv\!X^\dagger]$ such that 
$g_1' \stackrel{*}{\to}_{F} h$ and
$g_2' \stackrel{*}{\to}_{F} h$.\\

Suppose $g_1-g_2 \to_{f_i} g \stackrel{*}{\to}_{F} 0$ where
$g \stackrel{*}{\to}_{F} 0$ is a reduction sequence of length $n$.
\\

Let $t \in {\dashv\!X^\dagger}$ be the tagged term in $g_1-g_2$ to
which the reduction by $f_i$ is applied.
Let $u \in \dashv\!X^* \cup \{i\!d\}$, $v \in X^*$ such that 
$t=u \mathtt{LT}(f_i) v$, and let $k_1,k_2$ be the coefficients
of $t$ in $g_1,g_2$ respectively.
Now $k_1-k_2 \not= 0$ since it is the coefficient of $t$ in $g_1-g_2$.\\

Depending on whether $k_1$ and $k_2$ are zero or not we have the
following
zero- or one-step reductions:
$$g_1 \stackrel{=}{\to}_{f_i} g_1-k_1uf_iv, \quad
g_2 \stackrel{=}{\to}_{f_i} g_2-k_2uf_iv.$$

Since $g=(g_1-k_1uf_iv)-(g_2-k_2uf_iv)$ and $g \stackrel{*}{\to} 0$ in
$n$ 
steps, by the induction hypothesis there exists $h \in
K[\dashv\!X^\dagger]$
such that $g_1 - k_1uf_iv \stackrel{*}{\to}_{F} h$ and 
$g_2 - k_2uf_iv \stackrel{*}{\to}_{F} h$. Hence $g_1 \stackrel{*}{\to}_F
h$
and $g_2 \stackrel{*}{\to}_F h$.
\end{proof}

\begin{thm}[Test for confluence]
The reduction relation $\to_{F}$ generated by $F$ 
is complete on 
$K[\dashv \! X^\dagger]$ if and only
if all matches of $F$ resolve.
\end{thm}

\begin{proof}
In Lemma 4.2 we proved that $\to_{F}$ is Noetherian and therefore,
by Newman's Lemma for reduction relations on sets, we need only to prove
that
$\to_{F}$ is locally confluent.\\
  
  Let $f,g_1,g_2 \in K[\dashv\!X^\dagger]$ such that 
$f \to_{F} g_1$ and 
$f \to_{F} g_2$.
  Then 
$g_1=f-k_1u_1f_1v_1$ and
$g_2=f-k_2u_2f_2v_2$ for some $f_1,f_2 \in F$, $k_1,k_2 \in K$,
$u_1,u_2,v_1,v_2 \in \dashv\!X^* \cup X^*$.
  Let $m_1:=\mathtt{LT}(f_1)$ and $m_2:=\mathtt{LT}(f_2)$.\\
  
  If the reductions do not overlap on $f$, i.e. 
$u_1m_1v_1 \not= u_2m_2v_2$ then it is immediate that 
$g_1 \to_{F} h$ and 
$g_2 \to_{F} h$ where 
$h = f - k_1u_1f_1v_1 - k_2u_2f_2v_2$.\\

Otherwise $u_1m_1v_1=u_2m_2v_2$. 
In this case $m_1$ and $m_2$ may or may not coincide.\\

If they do not coincide, i.e. if there exists $w \in X^*$ such that
$u_1m_1v_1=u_1m_1wm_2v_2$ or $u_1m_1v_1=u_2m_2wm_1v_1$ then again
$g_1 \stackrel{*}{\to}_{F} h$ and 
$g_2 \stackrel{*}{\to}_{F} h$ where 
$h = f - u_1f_1wm_2v_2 - u_1m_1wf_2v_2$ or 
$h = f - u_2m_2wf_1v_1 - u_2f_2wm_1v_1$ respectively.\\

If the leading terms $m_1$ and $m_2$ do coincide then 
$u_1m_1v_1=u_2m_2v_2$ represents a multiple of a match between $f_1$ and
$f_2$,
i.e. there exist $u_1',u_2',v_1',v_2',w,z \in \dashv\!X^* \cup X^*$, 
such that $u_1=wu_1'$, $v_1=v_1'z$, $u_2=wu_2'$, $v_2=v_2'z$ and
$u_1'm_1v_1'=u_2'm_2v_2'$ represents a match between $f_1$ and $f_2$.
In this case $u_1'f_1v_1'-u_2'f_2v_2' \stackrel{*}{\to}_{F} 0$ by
assumption, and therefore 
$wu_1'f_1v_1'z-wu_2'f_2v_2'z 
  = u_1m_1v_1-u_2m_2v_2 \stackrel{*}{\to}_{F} 0$. By Lemma 4.5 this
  implies that there exists $h \in K[\dashv \! X^\dagger]$ such that
$g_1 \stackrel{*}{\to}_{F} h$ and 
$g_2 \stackrel{*}{\to}_{F} h$.\\

The converse of the above is easily checked. 
Suppose that $\to_{F}$ is confluent. Then any
S-polynomial arising from a match between polynomials is the result of
reducing one term in two different ways, i.e. 
$f \to_{F} g_1$ and
$f \to_{F} g_2$ for some $f,g_1,g_2\in K[X^\dagger]$. 
The S-polynomial is
equal to
$g_1- g_2$. The relation $\to_{F}$ is locally confluent and so
there
exists $h \in K[\dashv \! X^\dagger]$ such that 
$g_1 \to_{F} h$ and 
$g_2 \to_{F} h$. Therefore
$g_1 - g_2 \stackrel{*}{\to} h - h = 0$ as required.
\end{proof}

We may now apply the noncommutative version of Buchberger's algorithm 
(as described in \cite{TMora2}) to 
attempt to complete a mixed set of tagged and non-tagged polynomials.
To verify steps 4 and 5 of the algorithm we observe the following two 
technical lemmas.

\begin{lem}[Addition of S-polynomials]
Let $F=(F_T,F_P)$. If $f$ is an S-polynomial resulting from a match of $F$, 
then the congruences $\stackrel{*}{\leftrightarrow}_F$ and
$\stackrel{*}{\leftrightarrow}_{F\cup\{f\}}$ coincide.
\end{lem}
\begin{proof}
The result is proved by showing that, in each of the five cases, an
S-polynomial $f$ resulting from a match of polynomials $f_1,f_2 \in F$ 
can be written in the form $u_1f_1v_1-u_1f_2v_2$ and therefore 
$f \stackrel{*}{\leftrightarrow}_F 0$.
\end{proof}

\begin{lem}[Elimination of redundancies]
Let $F=(F_T,F_P)$. If $f\in F$ is such that $f \to_{F \setminus\{f\}} 0$ then 
the relations $\stackrel{*}{\to}_F$ and
$\stackrel{*}{\to}_{F\setminus\{f\}}$ coincide.
\end{lem}
\begin{proof}
The result is immediate, since for all $g \to_{\{f\}} h$ then 
$g = h+ kufv \stackrel{*}{\to}_{F\setminus\{f\}} h$ where $k \in K$, 
$ufv \in \dashv\!X^\dagger$.
\end{proof} 

\begin{alg}[Noncommutative Buchberger Algorithm with tags]

Given a set of tagged and non-tagged polynomials the algorithm 
attempts to complete the set with respect to a given ordering
so that the reduction relation generated is complete.
\begin{enumerate}
\item
(Input:) A mixed set of tagged and non-tagged polynomials 
$F=(F_T,F_P)$ where $F_T \subseteq K[\dashv\!X^\dagger]$ and 
$F_P \subseteq K[X^\dagger]$.
\item
(Initialise:)
Put $\mathtt{OLD} := F$ and $\mathtt{SPOL}:=\emptyset$. 
\item
(Search for matches:)
If the leading monomials of any of the polynomials overlap then calculate
the
resulting S-polynomial and attempt to reduce it using $\to_F$.
Record all non-zero reduced S-polynomials in the list $\mathtt{SPOL}$.
\item
(Add unresolved S-polynomials:)
When all matches have been considered define $\mathtt{NEW:=OLD \cup
SPOL}$.
\item
(Eliminate redundancies:)
Pass through $\mathtt{NEW}$ removing each polynomial in turn and
reducing it with respect to the other polynomials in $\mathtt{NEW}$. 
If a polynomial reduces to zero, delete it from
$\mathtt{NEW}$. Otherwise replace each with its reduced form.
\item
(Loop:)
Whilst $\mathtt{OLD} \not= \mathtt{NEW}$ set $\mathtt{OLD:=NEW}$,
$\mathtt{SPOL}:=\emptyset$ and return to step 3.
\item
(Output:)
A set $F:=\mathtt{NEW}$ of polynomials such that $\to_F$
is a complete reduction relation on $K[\dashv\!X^\dagger]$.
\end{enumerate}
\end{alg}

\begin{rem}[Left Ideals]
\emph{Placing tags to the right of polynomials rather than the left, 
i.e. working in $K[X^\dagger\vdash]$ or $K[X^*\vdash]$, by similar arguments
we can compute left ideals. The tags act as a block to multiplication: 
to calculate left ideals, one blocks the right multiplication with a tag on 
the right; to calculate right ideals, one blocks the left multiplication 
with a tag on the left. It is natural that two-sided ideals have no tags, 
since both multiplications are defined.}
\end{rem}

\begin{rem}[Implementation]
\emph{The use of the free monoid $(X \cup \{ \dashv \})^*$ is possible
in Definition 3.2, i.e. $f \to_F f-kuf_iv$ if 
$u\mathtt{LT}(f_i)v$ occurs in $f$ with coefficient $k$
$u,v \in (X \cup \{ \dashv \})^*$.  
If a match of any of the five types described above occurs between 
$f_1$ and $f_2$ then there exist
$u_1,u_2,v_1,v_2 \in (X \cup \{ \dashv \})^*$ such that
$u_1 m_1 v_1 =u_2 m_2 v_2$ (the converse is not true).
Unmeaningful monomials 
such as $\dashv\!xx\!\dashv\dashv\!x$ do not arise as a result of
any procedure we describe, including Algorithm 4.9.
Therefore there is no problem with the
computations taking place inside the free $K$-algebra 
$K[(X \cup\{\dashv\})^*]$.
This is useful computationally as it allows us to use a standard 
noncommutative Gr\"obner basis program as an implementation of the procedures.
In other words, this widens the scope of a noncommutative Gr\"obner basis 
program without modifying it: the program can now attempt to compute bases 
for one-sided ideals in finitely presented algebras.}
\end{rem}

\section{Application to Green's Relations}

The standard way of expressing the structure of an (abstract) semigroup
is in terms of Green's relations. The relations enable the expression
of the local structure of the semigroup in terms of groups with certain
actions 
on them. Eggbox diagrams depict the partitions of a semigroup into their
$L$-classes $R$-classes, $D$-classes and $H$-classes as defined by
Green's
relations. 
We can sometimes determine the classes by using Gr\"obner bases applied 
directly to the presentation. The
examples show that there is also the possibility of dealing with
infinite
semigroups having infinitely many $H$-classes, $L$-classes or
$R$-classes.
First we recall some definitions \cite{Howie}.\\

A nonempty subset $A$ of a
semigroup $S$ is a {\em right ideal} of $S$ if $AS \subseteq A$, where
$AS:=\{as:a \in A, s \in S\}$.
It is a {\em left ideal} of $S$ if $SA \subseteq A$.
If $x$ is an element of $S$ then the smallest right ideal of $S$
containing $x$
is $xS \cup \{x\}$, we denote this $\langle x\rangle^r $ as it is called
the
{\em right ideal generated by $x$}. Similarly the 
{\em left ideal generated by $x$} is $Sx \cup \{x\}$ and is denoted 
$\langle x \rangle^l$.\\

\textbf{Green's Relations}\\
Let $S$ be a semigroup and let $s$ and $t$ be elements of $S$.
We say that $s$ and $t$ are {\em L-related} if the left ideal
generated
by $s$ in $S$ is equal to the left ideal generated by $t$:
$$s \sim_L t \Leftrightarrow \langle s\rangle^l =\langle t\rangle^l .$$
Similarly they are {\em R-related} if the right ideals are the same:
$$s \sim_R t \Leftrightarrow \langle s\rangle^r =\langle t\rangle^r.$$

The $L$-relation is a right congruence on $S$ and the $R$-relation is a
left
congruence on $S$. (The right action of $S$ on itself is preserved by
the
mapping to the
$L$-classes - so $[x^y]_{\sim_L}=[xy]_{\sim_L}={[x]^y}_{\sim_L}$,
similarly for
the left action and $R$-classes.)
The elements $s$ and $t$ are said to be {\em H-related} if they are
\emph{both} $L$-related \emph{and} $R$-related, and are
{\em D-related} if they are \emph{either} $L$-related \emph{or}
$R$-related.\\

To determine whether $s$ and $t$ are $R$ (or $L$)-related
we can compute the appropriate Gr\"obner bases and compare them.
First let $K$ be (any) field.
Let $S$ have presentation $sgp\langle X|Rel\rangle$
Let $P$ be a Gr\"obner basis for $K[S]$ (so $K[X^\dagger]/\! =_P\cong
K[S]$).
We would add
the polynomial $\dashv \! s$ to the Gr\"obner basis system for $K[S]$
and
compute
the Gr\"obner basis, and see whether this was equivalent to the basis
obtained
for $\dashv \! t$.

\section{Examples}

Throughout the examples we will use the field $\mathbb{Q}$ and the 
standard length-lexicographical ordering $>$.

\begin{example}
\emph{The first example is a two element semigroup with presentation}
$sgp \langle x |x^3=x^2 \rangle.$\\

\emph{The Gr\"obner basis for the right ideal $\langle x\rangle^r$ is
$\{\dashv\!x, x^3-x^2\}$ and the Gr\"obner basis for $\langle x^2 \rangle^r$
is $\{\dashv\!x^2, x^3-x^2\}$. The Gr\"obner bases are different and therefore 
$x$ and $x^2$ are not $R$-related.
Similarly, the Gr\"obner basis for the left ideal $\langle x\rangle^l$  is
$\{x\!\vdash, x^3-x^2\}$ and the Gr\"obner basis for $\langle x^2\rangle^l$ is
$\{x^2\!\vdash, x^3-x^2\}$ so the elements are not $L$-related. Therefore
this semigroup has two $H$-classes.} 
\end{example}

\begin{example}
\emph{The following example is for the 
finite monoid $Sym(2)$ with semigroup 
presentation}

\vspace{-.4cm}
$$
sgp\langle e,s| e^2=e, s^3=s, s^2e=e, es^2=e, sese=ese, eses=ese\rangle.
$$ 
\vspace{-.5cm}

\emph{The Gr\"obner basis equivalent to the rewrite system is}

\vspace{-.4cm}
$$
F:=\{e^2-e, \ s^3-s, \ s^2e-e, \ es^2-e, \ eses-ese, \ sese-ese\}.
$$
\vspace{-.5cm}

\emph{The elements are $\{ e, s, es, se, s^2, ese, ses\}$, 
where $s^2$ is the identity element.
We calculate Gr\"obner bases for 
the right and left ideals for each of
the elements. The results are displayed 
in the table below. In detail, a Gr\"obner
basis for $\langle ses \rangle^r$ in $K[S]$ in 
$K[\dashv\! X^\dagger]$ is calculated by 
adding $\dashv\! ses$ to the set of polynomials 
$F$.
A match $s$ occurs on $\dashv\! sesse$ 
between $sse-e$ and $\dashv\! ses$. 
This results in the S-polynomial
$\dashv \! se(e)-(0)se$ 
which reduces to $\dashv \! se$. 
Another match of $es$ occurs on $\dashv \! seses$ 
between $eses-ese$ and $\dashv \! ses$. 
This results in the 
S-polynomial $\dashv \! s(ese)-(0)es$ 
which reduces to $\dashv \! ese$.
All further matches result in 
S-polynomials which reduce to zero. 
The polynomials we add to $F$ 
to obtain a Gr\"obner basis are 
$\{ \dashv \! se, \dashv \! ese \}$ 
(note that $\dashv \! ses$ is a 
multiple of $\dashv \! se$ so it is 
not required in the Gr\"obner basis).
The table lists the polynomials which, 
together with $F$, will give the 
Gr\"obner bases for the right and left 
ideals generated by single elements.}

\begin{center}
\begin{tabular}{l|l|l}
element & right ideal & left ideal\\
\hline

$e$    & $\dashv\! e                $    & $e \!\vdash
$\\
$s$    & $\dashv\! e,  \dashv\! s   $    & $e \!\vdash,   s \!\vdash
$\\
$es$   & $\dashv\! e                $    & $es \!\vdash , ese \!\vdash
$\\
$se$   & $\dashv\! se, \dashv\! ese $    & $e \!\vdash
$\\
$ss$   & $\dashv\! e,  \dashv\! s   $    & $e \!\vdash,   s \!\vdash
$\\
$ese$  & $\dashv\! ese              $    & $ese \!\vdash
$\\
$ses$  & $\dashv\! se, \dashv\! ese $    & $es \!\vdash , ese \!\vdash
$
\end{tabular}
\end{center}

\emph{Two elements whose right ideals are generated by the same
Gr\"obner basis have 
the same right ideal (similarly left), and so
it is immediately deducible that\\
the $R$-classes are $\{s,s^2\},\{e,es\},\{se,ses\}$ and $\{ese\}$, 
the $L$-classes are $\{s,s^2\},\{e,se\},\{es,ses\}$ and $\{ese\}$,
the $H$-classes are $\{s,s^2\},\{e\},\{se\},\{es\},\{ses\}$ and
$\{ese\}$ and
the $D$-classes are $\{s,s^2\},\{e,es,se,ses\}$ and $\{ese\}$.\\
The eggbox diagram is as follows
where $L$-classes are columns, $R$-classes are rows, 
$D$-classes are diagonal boxes
and $H$-classes are the small boxes:}
\\

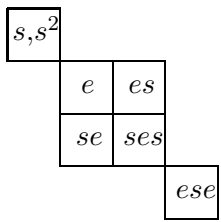
\setlength{\unitlength}{0.7cm}
\begin{picture}(10,4)
\put(9,3){\line(0,0){1}}
\put(10,1){\line(0,0){3}}
\put(11,1){\line(0,0){2}}
\put(12,0){\line(0,0){3}}
\put(13,0){\line(0,0){1}}
\put(9,4){\line(1,0){1}}
\put(9,3){\line(1,0){3}}
\put(10,2){\line(1,0){2}}
\put(10,1){\line(1,0){3}}
\put(12,0){\line(1,0){1}}
\put(9.1,3.4){$s,\!s^2$}
\put(10.3,1.4){$se$}
\put(10.4,2.4){$e$}
\put(11.2,1.4){$ses$}
\put(11.3,2.4){$es$}
\put(12.2,0.4){$ese$}
\end{picture}

\emph{This example could have been calculated by enumerating the
elements of each of the fourteen ideals -- a time consuming procedure which
calculates details which we do not require.}
\end{example}

\begin{example}
\emph{The next example is the Bicyclic monoid which is infinite.
We use the semigroup presentation} \;
$sgp \langle p,q, i | p i=p, q i= q, i p=p, iq=q, pq=i \rangle.$\\

\emph{The equivalent Gr\"obner basis, defined on $K[\{p,q,i\}^\dagger]$, is 
$\{p i-p, q i- q, i p-p, iq-q, pq-i\}$.
We begin the table as before:}

\begin{center}
\begin{tabular}{l|l|l}
element & right ideal & left ideal\\
\hline
$i\!d$    & $\dashv\! i$. & $i\! \vdash $.\\
$p$      & $\dashv\! i$. & $p\! \vdash $.\\
$q$      & $\dashv\! q$.   & $q\! \vdash $.\\
$p^2$    & $\dashv\! i$. & $p^2\! \vdash $.\\
$qp$     & $\dashv\! q$.   & $p\! \vdash $.\\
$q^2$    & $\dashv\! q^2$. & $i\! \vdash $.\\
$\cdots$ & $\cdots$        & $\cdots$\\
$q^np^m$ & $\dashv\! q^n$. & $p^m\! \vdash $.\\
\end{tabular}
\end{center}
\emph{It can be seen that there are infinitely many $L$-classes 
and infinitely many $R$-classes.
Representatives for the $L$-classes are the elements of $\{q\}^*$
because 
$q^np^m\! \vdash  \to q^n\! \vdash $ (using the S-polynomial resulting
from
$p^n(q^np^m\! \vdash )\to p^n\! \vdash $ with  $(p^nq^n)p^m\! \vdash
\to p^m\! \vdash $).
Similarly the elements of $\{p\}^*$ are representatives for the
$R$-classes.
All elements are $D$-related and none of them are $H$-related. So the
eggbox
diagram would be an infinitely large box of cells with one element in
each
cell. This means that the monoid is} bisimple.
\end{example}

\begin{example}
\emph{Now consider the Polycyclic monoid $P_n$ which has monoid
presentation}

\vspace{-.4cm}
$$
mon\langle x_1,\ldots,x_n,y_1,\ldots,y_n,o, i\!d
\; | \; ox_i \! = \! x_io \! = \! oy_i \! = \! y_io \! = \! o, x_iy_i \!
= \! i\!d,
x_iy_j \! = \! o  \text{\emph{ for }} i,j \! = \! 1,\ldots,{n\!-\!1},
i\not \! = \! j\rangle
$$
\vspace{-.4cm}

\emph{and therefore the Gr\"obner basis for $K[P_n]$, where $K$ is a
field, is}

\vspace{-.2cm}
$$
\{x_iy_i-i\!d, x_iy_j-0 \text{ for } i,j=1,\ldots,{n-1},\, i \not= j \}.
$$
\vspace{-.2cm}

\emph{Green's relations for the polycyclic monoids are naturally similar to
those for the Bicyclic monoid. The $L$-classes are represented by
sequences of
$y_i$'s and the
$R$-classes are represented by sequences of $x_i$'s. To verify this, let
$X=x_{i_1}\cdots x_{i_k}$ be a
general word in the $x_i$'s, and let $Y$ be $y_{j_1}\cdots y_{j_l}$ a
general
word in the $y_j$'s.
Then we can show that $YX \sim_L X$.
To do this consider $\langle YX\! \vdash  \rangle$ and $\langle X\!
\vdash  \rangle$.
 To find a Gr\"obner basis for $\langle YX\! \vdash  \rangle$ consider
the match
 $x_{j_l} \cdots x_{j_1}y_{j_1} \cdots y_{j_l}x_{i_1} \cdots
x_{i_k}\! \vdash $.
This results in the S-polynomial
 $(i\!d)x_{i_1} \cdots x_{i_k}\! \vdash  - x_{j_l} \cdots x_{j_1}(0)$
which simplifies to $x_{i_1} \cdots x_{i_k}\! \vdash  = X\! \vdash $.
This is
a
Gr\"obner basis for $\langle YX\! \vdash  \rangle$, and so $\langle YX\!
\vdash  \rangle
=\langle
X\! \vdash  \rangle$. Similarly $\langle \dashv \! YX \rangle= \langle
\dashv \! Y \rangle$ so
$YX
\sim_R X$ for any $Y = y_{j_1}\cdots y_{j_l}$, $X = x_{i_1}\cdots
x_{i_k}$.}
\\

\emph{The eggbox diagram is drawn below. As before the $L$ classes are the
columns
and the $R$-classes the rows, $H$-classes are the cells, and there is
just one
$D$-class other than the one containing the zero. This proves that the
polycyclic monoids are bisimple.
The diagram is more conventional than the
previous one, as classes are listed but not individual elements, instead
the
number of elements in each cell is indicated.}

\setlength{\unitlength}{0.8cm}
\begin{picture}(16,12)
\put(5,10){\line(0,0){1}}
\put(6,1){\line(0,0){1}}
\put(6,4){\line(0,0){2}}
\put(6,7){\line(0,0){4}}
\put(7,1){\line(0,0){1}}
\put(7,4){\line(0,0){2}}
\put(7,7){\line(0,0){3}}
\put(8,1){\line(0,0){1}}
\put(8,4){\line(0,0){2}}
\put(8,7){\line(0,0){3}}
\put(9,1){\line(0,0){1}}
\put(9,4){\line(0,0){2}}
\put(9,7){\line(0,0){3}}
\put(10,1){\line(0,0){1}}
\put(10,4){\line(0,0){2}}
\put(10,7){\line(0,0){3}}
\put(11,1){\line(0,0){1}}
\put(11,4){\line(0,0){2}}
\put(11,7){\line(0,0){3}}
\put(12,1){\line(0,0){1}}
\put(12,4){\line(0,0){2}}
\put(12,7){\line(0,0){3}}
\put(14,1){\line(0,0){1}}
\put(14,4){\line(0,0){2}}
\put(14,7){\line(0,0){3}}
\put(15,1){\line(0,0){1}}
\put(15,4){\line(0,0){2}}
\put(15,7){\line(0,0){3}}
\put(5,10){\line(1,0){4}}
\put(5,11){\line(1,0){1}}
\put(6,9){\line(1,0){3}}
\put(6,8){\line(1,0){3}}
\put(6,7){\line(1,0){3}}
\put(6,6){\line(1,0){3}}
\put(6,5){\line(1,0){3}}
\put(6,4){\line(1,0){3}}
\put(6,2){\line(1,0){3}}
\put(6,1){\line(1,0){3}}
\put(10,9){\line(1,0){2}}
\put(10,8){\line(1,0){2}}
\put(10,7){\line(1,0){2}}
\put(10,6){\line(1,0){2}}
\put(10,5){\line(1,0){2}}
\put(10,4){\line(1,0){2}}
\put(10,2){\line(1,0){2}}
\put(10,1){\line(1,0){2}}
\put(10,10){\line(1,0){2}}
\put(14,9){\line(1,0){1}}
\put(14,8){\line(1,0){1}}
\put(14,7){\line(1,0){1}}
\put(14,6){\line(1,0){1}}
\put(14,5){\line(1,0){1}}
\put(14,4){\line(1,0){1}}
\put(14,2){\line(1,0){1}}
\put(14,1){\line(1,0){1}}
\put(14,10){\line(1,0){1}}
\put(4.2,1.4){$[X]$}
\put(4.2,4.4){$[x_1x_2]$}
\put(4.2,5.4){$[{x_1}^2]$}
\put(4.2,7.4){$[x_2]$}
\put(4.2,8.4){$[x_1]$}
\put(4.2,9.4){$[i\!d]$}   \put(4.2,10.4){$[0]$} \put(5.2,11.4){$[0]$}
\put(6.2,11.4){$[i\!d]$}
\put(7.2,11.4){$[y_1]$}
\put(8.2,11.4){$[y_2]$}
\put(10.2,11.4){$[{y_1}^2]$}
\put(11.2,11.4){$[y_1y_2]$}
\put(14.2,11.4){$[Y]$}
\put(6,0.3){\line(0,0){0.4}}
\put(7,0.3){\line(0,0){0.4}}
\put(8,0.3){\line(0,0){0.4}}
\put(9,0.3){\line(0,0){0.4}}
\put(10,0.3){\line(0,0){0.4}}
\put(11,0.3){\line(0,0){0.4}}
\put(12,0.3){\line(0,0){0.4}}
\put(14,0.3){\line(0,0){0.4}}
\put(15,0.3){\line(0,0){0.4}}
\put(6,2.3){\line(0,0){0.4}}
\put(7,2.3){\line(0,0){0.4}}
\put(8,2.3){\line(0,0){0.4}}
\put(9,2.3){\line(0,0){0.4}}
\put(10,2.3){\line(0,0){0.4}}
\put(11,2.3){\line(0,0){0.4}}
\put(12,2.3){\line(0,0){0.4}}
\put(14,2.3){\line(0,0){0.4}}
\put(15,2.3){\line(0,0){0.4}}
\put(6,3.3){\line(0,0){0.4}}
\put(7,3.3){\line(0,0){0.4}}
\put(8,3.3){\line(0,0){0.4}}
\put(9,3.3){\line(0,0){0.4}}
\put(10,3.3){\line(0,0){0.4}}
\put(11,3.3){\line(0,0){0.4}}
\put(12,3.3){\line(0,0){0.4}}
\put(14,3.3){\line(0,0){0.4}}
\put(15,3.3){\line(0,0){0.4}}
\put(6,6.3){\line(0,0){0.4}}
\put(7,6.3){\line(0,0){0.4}}
\put(8,6.3){\line(0,0){0.4}}
\put(9,6.3){\line(0,0){0.4}}
\put(10,6.3){\line(0,0){0.4}}
\put(11,6.3){\line(0,0){0.4}}
\put(12,6.3){\line(0,0){0.4}}
\put(14,6.3){\line(0,0){0.4}}
\put(15,6.3){\line(0,0){0.4}}
\put(9.3,1){\line(1,0){0.4}}
\put(9.3,2){\line(1,0){0.4}}
\put(9.3,4){\line(1,0){0.4}}
\put(9.3,5){\line(1,0){0.4}}
\put(9.3,6){\line(1,0){0.4}}
\put(9.3,7){\line(1,0){0.4}}
\put(9.3,8){\line(1,0){0.4}}
\put(9.3,9){\line(1,0){0.4}}
\put(9.3,10){\line(1,0){0.4}}
\put(12.3,1){\line(1,0){0.4}}
\put(12.3,2){\line(1,0){0.4}}
\put(12.3,4){\line(1,0){0.4}}
\put(12.3,5){\line(1,0){0.4}}
\put(12.3,6){\line(1,0){0.4}}
\put(12.3,7){\line(1,0){0.4}}
\put(12.3,8){\line(1,0){0.4}}
\put(12.3,9){\line(1,0){0.4}}
\put(12.3,10){\line(1,0){0.4}}
\put(13.3,1){\line(1,0){0.4}}
\put(13.3,2){\line(1,0){0.4}}
\put(13.3,4){\line(1,0){0.4}}
\put(13.3,5){\line(1,0){0.4}}
\put(13.3,6){\line(1,0){0.4}}
\put(13.3,7){\line(1,0){0.4}}
\put(13.3,8){\line(1,0){0.4}}
\put(13.3,9){\line(1,0){0.4}}
\put(13.3,10){\line(1,0){0.4}}
\put(15.3,1){\line(1,0){0.4}}
\put(15.3,2){\line(1,0){0.4}}
\put(15.3,4){\line(1,0){0.4}}
\put(15.3,5){\line(1,0){0.4}}
\put(15.3,6){\line(1,0){0.4}}
\put(15.3,7){\line(1,0){0.4}}
\put(15.3,8){\line(1,0){0.4}}
\put(15.3,9){\line(1,0){0.4}}
\put(15.3,10){\line(1,0){0.4}}   \put(5.4,10.4){1}
\put(6.4,9.4){1} \put(6.4,8.4){1} \put(6.4,7.4){1}
\put(6.4,5.4){1} \put(6.4,4.4){1} \put(6.4,1.4){1}
\put(7.4,9.4){1} \put(7.4,8.4){1} \put(7.4,7.4){1}
\put(7.4,5.4){1} \put(7.4,4.4){1} \put(7.4,1.4){1}
\put(8.4,9.4){1} \put(8.4,8.4){1} \put(8.4,7.4){1}
\put(8.4,5.4){1} \put(8.4,4.4){1} \put(8.4,1.4){1}
\put(10.4,9.4){1} \put(10.4,8.4){1} \put(10.4,7.4){1}
\put(10.4,5.4){1} \put(10.4,4.4){1} \put(10.4,1.4){1}
\put(11.4,9.4){1} \put(11.4,8.4){1} \put(11.4,7.4){1}
\put(11.4,5.4){1} \put(11.4,4.4){1} \put(11.4,1.4){1}
\put(14.4,9.4){1} \put(14.4,8.4){1} \put(14.4,7.4){1}
\put(14.4,5.4){1} \put(14.4,4.4){1} \put(14.4,1.4){1}
\end{picture}
\end{example}

These examples illustrate the fact that Buchberger's algorithm can be used to
compute Green's relations for (infinite) semigroups which have finite
complete presentations. 
Previous methods for calculating minimal ideals from presentations of
semigroups were variations on the classical Todd-Coxeter enumeration
procedure \cite{Campbell95}.
This is an alternative computational approach 
to that given in \cite{Linton,Linton2} which uses the transformation 
representation of a semigroup rather than a presentation.
As with \cite{Linton2} the methods described in this paper provide for local computations, 
concerning a single $R$-class, without computing the whole semigroup.
The one-sided Gr\"obner basis methods have
limitations in that a complete rewrite system with respect to the chosen
order might not be found, but they do give the possibility of calculating the
structure of infinite semigroups and do not require the determination of a
transformation representation for those semigroups which arise 
naturally as presentations.\\

The calculations of the examples were achieved using a $\mathsf{GAP3}$ 
implementation of the Gr\"obner basis procedures  for polynomials in 
noncommutative variables over $\mathbb{Q}$ as described in \cite{TMora}. 
Further details of this program can be found in \cite{Annethesis} or e-mail
the author.
Other implementations (e.g. OPAL, Bergman) are more powerful: the key point 
of this paper is to point out that such programs can be used 
for a wider range of problems than has previously been recorded.\\


{\small
\begin{thebibliography}{99}
\bibitem{TAT} F.Baader and T.Nipkow : 
                Term Rewriting and All That,
                \emph{Cambridge University Press} 1998.

\bibitem{RISC} B.Buchberger and F.Winkler :  
                ``Gr\"obner Bases and Applications'',
                \emph{``33 Years of Gr\"obner Bases'' RISC-Linz 2-4 Feb 1998, Proc. London Math. Soc. vol.251} 1998.
 
\bibitem{Campbell95} C.M.Campbell, N.Ru\v skuc, E.F.Robertson and R.M.Thomas :  
               ``Rewriting a Semigroup Presentation'',
                \emph{International Journal of Algebra and Computation, 
                vol.5 no.1 p81-103} 1995.

\bibitem{Annethesis} A. Heyworth, 
     ``Rewriting and Noncommutative Gr\"obner Bases with Applications to
     Kan Extensions and Identities Among Relations'',
     \emph{UWB Math Preprint 98.23},
     {\tt http://xxx.soton.ac.uk/abs/math/9812097}, 1998.
     
\bibitem{paper3} A. Heyworth, 
     ``Rewriting as a Special Case of Noncommutative Gr\"obner Basis Theory'',
     \emph{UWB Math Preprint 98.22}, 
     {\tt http://xxx.soton.ac.uk/abs/math/9901044}, 1998. 
               
\bibitem{Howie} J.M.Howie:
                ``Fundamentals of Semigroup Theory'',
                \emph{LMS new series vol.12, Oxford Science Publications} 1995.
              
\bibitem{Linton} S.A.Linton, G.Peiffer, E.F.Robertson and N.Ru\v skuc:
                ``Groups and Actions in Transformation Semigroups,''
                \emph{Mathematische Zeitschrift vol.228 p435-450}, 1998.               
                                
\bibitem{Linton2} S.A.Linton, G.Peiffer, E.F.Robertson and N.Ru\v skuc :
                ``Computing Transformation Semigroups'',
                \emph{Journal of Symbolic Computation vol.11 p1-18}, 1998.                

\bibitem{FMora} F.Mora :
                ``Gr\"obner Bases for Noncommutative Polynomial Rings'',
                {\em Proc. AAECC-3, LNCS 229, p353-362, Springer}, 1985.
 
\bibitem{TMora} T.Mora : 
                ``Gr\"obner Bases and the Word Problem'',
                \emph{Preprint, University of Genova} 1987.

\bibitem{TMora2} T.Mora : 
                ``An Introduction to Commutative and Noncommutative
                Gr\"obner Bases'',
                \emph{Theoretical Computer Science vol.134 p131-173} 1994.

\bibitem{Birgit} B.Reinert : 
                ``On Gr\"obner Bases in Monoid and Group Rings'',
                \emph{PhD Thesis, Universit\"at Kaiserslautern} 1995.

\bibitem{MRC} B.Reinert and D.Zecker : 
                ``MRC - A System for Computing Gr\"obner Bases in
                Monoid and Group Rings'',
                \emph{Universit\"at Kaiserslautern Preprint} 1998.

\bibitem{Ufn} V.Ufnarovski :
                ``Introduction to Noncommutative Gr\"obner Basis Theory'',
                 \emph{in Gr\"obner Bases and Applications, 
                 B.Buchberger and F.Winkler (eds),
                 Proc. London Math. Soc. vol.251 p305-322} 1998.

\end{thebibliography}
}
\end{document}